\documentclass[12pt]{amsart} 
\usepackage{latexsym,fancyhdr,amssymb,color,amsmath,amsthm,graphicx,listings,comment}

\usepackage[section]{placeins}      \let\paragraph\subsection
\title{More on Numbers and Graphs}
\author{Oliver Knill} \date{5/30/2019}
\address{Department of Mathematics \\ Harvard University \\ Cambridge, MA, 02138 }

\begin{document}
\begin{abstract}
In this note we revisit a ``ring of graphs" 
$\mathcal{Q}$ in which the set of finite simple graphs $\mathcal{N}$ extend the role of the 
natural numbers $\mathbb{N}$ and the signed graphs $\mathcal{Z}$ extend the role of 
the integers $\mathbb{Z}$. We point out the existence of a norm which allows to complete 
$\mathcal{Q}$ to a real or complex Banach algebra $\mathcal{R}$ or $\mathcal{C}$.
\end{abstract}

\maketitle

\section{A normed ring in a nutshell}

\paragraph{}
The {\bf Zykov-Sabidussi ring} $(\mathcal{Z},+,\star)$ with 
{\bf Zykov join} $(V,E) + (W,F) = (V \cup W, E \cup F \cup \{ (a,b), a \in V, b \in W \}$ 
as addition and {\bf Sabidussi multiplication} 
$(V,E) \star (W,F) =(V \times W, \{ (a \times b,c \times d)$, $(a,c) \in E$ or $(b,d) \in F \}$ 
imposes an associative and commutative ring structure on the class $\mathcal{Z}$ of signed 
finite simple graphs $G=(V,E)$, where $V$ is the vertex set and $E$ the edge set. 
The distributivity $A*(B+C)=A*B+A*C$ can be easily checked. This goes back to work of 
Zykov \cite{Zykov} and Sabidussi \cite{Sabidussi}. See also \cite{HammackImrichKlavzar}.
The construction of $\mathcal{Z}$ is analogue to the construction of integers from natural numbers:
the monoid $(\mathcal{N},+)$ of finite simple graphs is first Grothendieck completed to a group $(\mathcal{Z},+)$,
then the multiplication $\star$ is extended from $(\mathcal{N},+,\star)$ to a ring $(\mathcal{Z},+,\star)$. 
The class of complete graphs generates then a subring of $\mathcal{Z}$ that is isomorphic to $\mathbb{Z}$ so that 
$\mathcal{Z}$ naturally extends integer arithmetic of $\mathbb{Z}$ to the larger ring $\mathcal{Z}$ of geometries. 

\begin{figure}[!htpb]
\scalebox{0.371}{\includegraphics{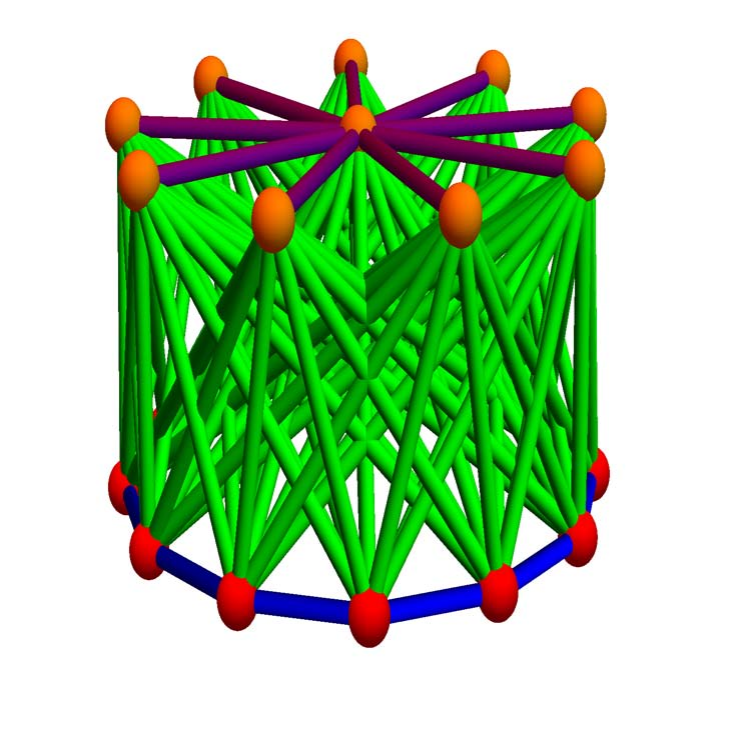}}
\scalebox{0.371}{\includegraphics{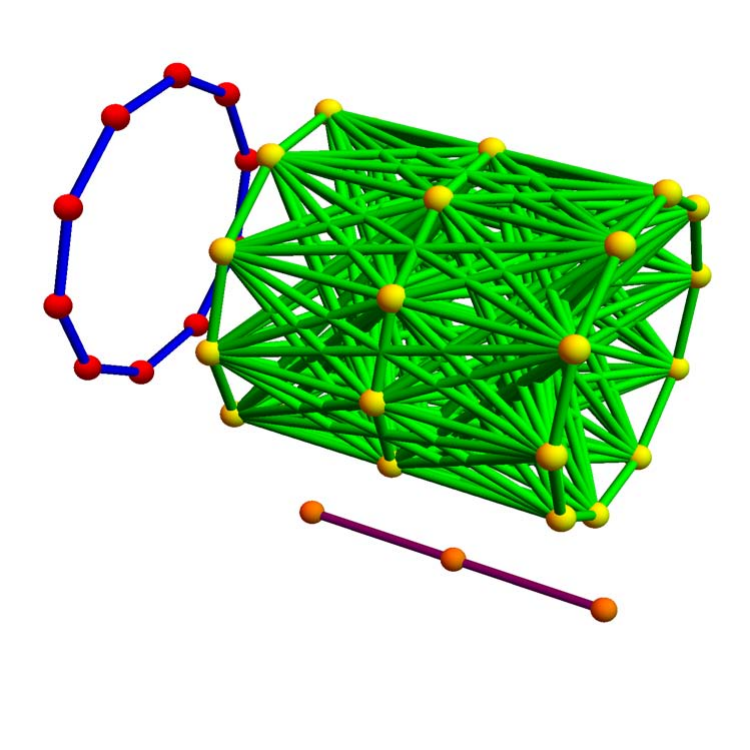}}
\label{addmultiply}
\caption{
Addition and multiplication of graphs. The sum is a graph on the vertex union,
the product is a graph on the vertex Cartesian product. The operations are compatible
with the clique number $c$ which is 1 plus the maximal dimension of the clique 
complex generated by the graph. We have $c(G*H)=c(G) c(H)$, $c(G+H)=c(G)+c(H)$.
}
\end{figure}

\paragraph{}
The {\bf clique number} $c(G)$ which satisfies $c(A+B) = c(A) + c(B), c(A*B)=c(A) c(A)$ for graphs $A,B \in \mathcal{N}$
extends then to the additive group $(\mathcal{Z},+)$ by setting $c(G)=c(A)-c(B)$ if $G=A-B$. We think of it as a trace. 
The function $c$ extends to the ring $(\mathcal{Z},+,\star)$
and so, by defining $c(A/B)=c(A)/c(B)$ also onto the ring $\mathbb{Q} = \{ A/B \; | \; A,B \in \mathcal{Z}, c(B) \neq 0 \}$. 
The relation $c(A+B)=c(A)+c(B), c(A*B)=c(A)c(B)$ continues to hold in the {\bf field} $\mathbb{Q}$ inside $\mathcal{Q}$ 
that is generated by complete graphs. The set $\mathcal{Q}$ is a ring because if two fractions $A/B$ and $C/D$ are in 
$\mathcal{Q}$, then the product and sum is in $\mathcal{Q}$ again is there because $c(BD)= c(B) c(D) \neq 0$. 
The subset $\mathbb{Q}$ generated by complete graphs is even a field as for any non-zero element, we can find an inverse.
Its completion $\mathbb{R}$ of $\mathbb{Q}$ then produce {\bf scalars} which render $\mathcal{Q}$ to a vector space. 

\paragraph{}
Define on $\mathcal{Z}$ the {\bf norm} $|G|=\inf_{A-B \sim G} \; c(A)+c(B)$
on $\mathcal{Z}$, where $A,B$ are in $\mathcal{N}$. 
If $G=A-B$ is a representation for which $A,B$ have no common additive factor, the
norm is just $c(A)+c(B)$. (This is similar as if we would define a degree of a rational function as the 
sum of the degrees of nominator and denominator. This means that we write the rational function in a reduced
form and define the degree as the sum of that reduced fraction.)
The norm extends with $|A/B| = |A|/|B|$ from $\mathcal{Z}$ to $\mathcal{Q}$. 
It induces the usual absolute value on $\mathbb{Q}$. 
Now, unlike $c(G)$ which can be zero also for non-zero $G$, the
norm $|G|$ is zero only if $G$ is zero. 

\paragraph{}
Because $|A \star B| \leq |A| |B|$, topologically completing the ring $\mathcal{Q}$ with 
respect to the metric defined by $|\cdot|$ 
gives now a {\bf commutative Banach algebra} $\mathcal{R}$. 
The completion $\mathcal{Q} \to \mathcal{R}$ is analogue to
completing the rationals $\mathbb{Q}$ to the reals $\mathbb{R}$. Now, the finite simple graphs $\mathcal{N}$
play the role of $\mathbb{N} \cup \{0\}$ but $\mathcal{Q}$ is not a field.
The arithmetic $(\mathcal{R},+,*)$ and the norm $| \cdot |$
can be extended to $\mathcal{C} = \mathcal{R} + i \mathcal{R}$, where it gives a {\bf complex Banach algebra}.
The units in $\mathcal{Q}$ are the non-empty signed graphs without edges, so that the multiplicative group
of units in $\mathcal{Q}$ is isomorphic to $(\mathbb{Z} \setminus \{0\},*)$ in standard arithmetic.

\paragraph{}
We are not aware of an other ring structure on signed graphs which leads to a Banach algebra.
In order to have {\bf real or complex vector space structure} on a space of graphs we need an action
$A \to \lambda A$ of the reals on the space of graphs. In our case, this is already given as the reals can be seen as
part of the completion of the ring. It is a bonus that the vector space is an algebra and even a Banach algebra.
It is possible because the product $*$ is compatible with the clique number functional $c$
which is used to define the norm. 
This means that the scalar multiplication $A \to \lambda A$ works and satisfies $|\lambda A| = |\lambda| |A|$.

\paragraph{}
The $f$-function $f_G(t)$ and Euler characteristic can be extended to $\mathcal{Z}$ (allowing also the value infinity) by
$f_{A-B}(t)=f_A(t)/f_B(t)$ and $\chi(A-B)=1-f_{A-B}$. As $f_B(-1)$ can be zero, it is better to keep the rational 
function instead).
Things do not go further, as $\chi(A*B)$ can not be expressed through $f_A$ and $f_B$.

\paragraph{}
It appears that we need the strong ring \cite{StrongRing} with disjoint union and Cartesian product
to have multiplicative compatibility with Euler characteristic. In the dual addition which is the disjoint union,
the Euler characteristic is obviously additive. In the addition $+$ given by the join, the genus $g(G) = 1-\chi(G)$ is
multiplicative $g(G+H)=g(G) g(H)$ as the $f$-function is multiplicative $f_{G+H}(t) = f_G(t) f_H(t)$.
But the strong ring leaves the class of graphs and is homotopy equivalent to a product given by 
the {\bf Stanley-Reisner ring} in which every graph is represented as a polynomial in 
$f_0=|V(G)|$ variables of degree $c(d)$. The graph $K_3$ for
example is $xyz+xy+yz+xz+x+y+z$. Addition and multiplication in this larger ring however produces 
objects which are no more graphs. We stayed within graphs by defining the product differently as in
\cite{KnillKuenneth}. 

\section{The clique number and arithmetic}

\paragraph{}
Among the many binary operations available on graphs, the join $+$ operation is a particular nice one 
\cite{Zykov}. It is the analog of the join in topology and interesting for various reasons: 
first of all, it imposes a monoid structure on finite simple graphs 
which has as a sub-monoid the set of {\bf spheres} and the monoid of {\bf Dehn-Sommerville graphs} (defined 
recursively as $X_{-1}=\{ 0\}$ and $X_{d+1} = \{ G \in \mathcal{N}, \chi(G)=1+(-1)^d$ and $S(x) \in X_{d}$ for all vertices $x \in V(G) \}$.
As in general, $1-\chi(G+H) = (1-\chi(G)) (1-\chi(Y))$ and 
$S_{G+H}(x) = S_G(x) + H$ for $x \in G$ and $S_{G+H}(y) = G+S_H(y)$ for $y \in H$, the class $X$ of
Dehn-Sommerville graphs is closed under addition, a fact which is useful in combinatorial topology.)
Especially adding the zero-dimensional sphere $S^0$, the {\bf suspension} $G \to G+S_0$ 
is an important operation. 

\paragraph{}
The join is also an operation which is compatible with the {\bf clique number} $c(G)$ (counting the cardinality of 
the vertex set of the largest complete subgraph in $G$), in the sense that 
$$ c(A+B)=c(A)+c(B) $$ 
for any two graphs $A,B$. The join operation is also interesting as it is dual to the {\bf disjoint union} operation 
of graphs. The {\bf graph complement} of a graph $G=(V,E)$ is $(V,E')$, where $E \cup E'$ is the edge set of the 
complete graph on $V$ and $E \cap E' = \emptyset$. Now, $\overline{G + H} = \overline{G} \cup \overline{H}$, 
where $\cup$ is the disjoint union. It follows that the join monoid is an additive {\bf unique factorization domain}
in which the {\bf additive primes} are the graphs $G$ for which $\overline{G}$ is connected. The prime
factorization of a graph $G$ is then obtained by identifying the connected components of $\overline{G}$. One has to
get a bit used to this, as in the classical natural numbers $\mathbb{N}$, there is only one additive prime $1$, 
because any natural number can be written uniquely as a sum of $1$. 
For graphs $G \in \mathcal{N}$, there are many additive primes: 
any graph for which the graph complement is connected is an additive prime.  In $\mathbb{N}$ of course there is
only one additive prime: $1$. Additive counting is the story of $1$. 

\paragraph{}
Any associative multiplication on graphs compatible with the disjoint union in the dual picture 
defines a multiplication compatible with the join. 
There are many associative ones \cite{HammackImrichKlavzar}: the weak product, the tensor product and 
the strong product are just three of them. 
But only the strong product and its dual has the property that it is compatible with the 
clique number $c(G)$ giving the cardinality of the vertex set of the largest complete subgraph 
of $G$. The disjoint union itself is not compatible as it satisfies 
$c(G \cup H) = {\rm max}(c(G),c(H))$ for the disjoint union $\cup$.

\paragraph{}
The ring $(\mathcal{Z},+,\star)$ naturally extends {\bf integer arithmetic} to a larger 
geometric frame work as also the multiplication is compatible with dimension: 
$$  c(A * B) = c(A) c(B) \; . $$ 
In \cite{ArithmeticGraphs} we suggested to look at the {\bf field of fractions} containing the ring $\mathcal{Z}$
(which is possible because $\mathcal{Z}$ is an integral domain, actually every additive prime in the ring has a unique
multiplicative prime factorization). Counter examples to unique prime factorization appeared in 
\cite{ImrichKlavzar,HammackImrichKlavzar} (and \cite{StrongRing} in the dual case) and is related on the fact that 
$\mathbb{N}[x]$ has no unique prime factorization: $(1+x+x^2)(1+x^3) =(1+x^2+x^4)(1+x)$.

\paragraph{} 
We prefer to to not take the field of fractions as we can not maintain the field
property in the completion anyway. Here, we define the quotient $A/B$ as long as $c(B)$ is not zero. 
The ring $\mathcal{Z}$ contains the subring $\mathcal{K}$ of signed complete graphs $K_n$ 
enhanced with $0$, which is the empty graph. The subring in the usual way completes to a 
field, which is is isomorphic to the real numbers $\mathbb{R}$. The construction of that field is identical 
to the construction for the usual arithmetic. 

\paragraph{}
Define the {\bf norm}
$$   |G|={\rm inf}_{A-B \sim G} \; c(A)+c(B)  \; ,  $$ 
if $A,B$ are finite simple graphs
(meaning that there exists no $C$ with $A=A'+C$ and $B=B'+C$). 
The definition might look a bit strange at first as this is not necessary for
integers where we can always write either $n$ or $-n$ for an integer and define $|\pm n|=n$. 
In the ring of signed graphs, it is possible that $A-B$ can not be simplified. It can only be
simplified if $B$ is an additive factor of $A$ or $A$ is an additive factor of $B$. 

\paragraph{}
This translation invariant norm leads to a metric $d(G,H)=|G-H|$. The
{\bf triangle inequality} can be verified readily: if $G=A-B$ and $H=C-D$ for $A,B,C,D \in \mathcal{N}$,
then $|G|=|A|+|B|$ and $|H|=|C|+|D|$ but $G-H=(A+C)-(B+D)$ might have a common additive factor
$U$ such that $G-H=(A+C-U) - (B+D-U)$ has a smaller norm. 
The triangle inequality is in general strict: for $G=C_4-C_5$ and $H=C_5-C_6$ for example,
we have (because $c(C_n)=2$ as a complex with maximal dimension $1$), 
$|G|=|H|=4$ and $|G-H|=|C_4-C_6|=4$ which is strictly smaller than $|G|+|H|=8$. 
For graphs, the triangle inequality $|A-B| + |B-C| \geq |A-C|$ 
can be an exact equality. This happens for example if $B=0$ and $A,C$ do not have a common additive factor. 
We also have $|F-G|=|G-F|$ and $|F-G|=0$ if and only if $c(F')+c(G')=0$ for some compatible
$F'=F+U,G'=G+U$ with $F'=G'=0$ which is the case if and only if $F=G$.

\paragraph{}
As $|A*B| \leq |A| |B|$, the topological ring $(\mathcal{Q},+,*,|\cdot|)$ 
can be completed to become a Banach space which is then an 
infinite-dimensional Banach algebra $(\mathcal{R},+,*,|\cdot|)$. 
First of all it is a linear space. There is an addition on $\mathcal{Q}$ which extends to $\mathcal{R}$.
Then as there is an incarnation of $\mathbb{Q}$ or $\mathbb{R}$
inside the class of graphs generated by complete graphs we have a {\bf scalar multiplication},
which is needed to define a linear space. The norm $|\cdot|$ makes it a normed linear space. 

\begin{comment}
\paragraph{}
There is no hope for $|A*B| = |A| \,  |B|$ in general, as this would lead to an 
associative and commutative {\bf division algebra}, a structure which can only be either 
$\mathbb{R}$ or $\mathbb{C}$. 
To get elements $A,B$ with $A*B=0$ one can look for $A=U-V, B=U+V$ for two $U,V$ satisfying 
$U^2=V^2$ which are not $U=\pm V$. 
\end{comment}

\section{Euler characteristic}

\paragraph{}
The {\bf Euler characteristic} for a finite simple graph $G \in \mathbb{N}$ 
is defined as $1-f_G(-1)=f_0-f_1+f_2- \cdots = 1-f_G(-1)$, where $f_G$ is the 
{\bf $f$-function} $f(t) = 1+f_0 t + f_1 t^2 + \cdots$ of $G$ defined by the 
$f$-vector $f=(f_0,f_1, \dots)$. The $f$-function can be extended to the full group
$(\mathcal{Z},+)$ of signed graphs by defining
$$  f_{A-B}(t)= \frac{f_A(t)}{f_B(t)} \;  $$
which is a rational function. One could define Euler characteristic on $\mathcal{Z}$
by putting 
$\chi(A-B) = 1- \frac{f_A(-1)}{f_B(-1)}$ if $f_B(-1)$ is not zero. Since $f_B(-1)=0$ 
is possible (it happens whenever the Euler characteristic of $B$ is $1$), 
it might be better just to extend the $f$ function to $\mathcal{Z}$ and not $\chi(G)$. 
Is there an extension to the ring $(\mathcal{Q},+,\star)$? Let's see: 

\paragraph{}
On the field generated by complete graphs, it is possible to extend Euler characteristic to the entire
field. But it does not lead to anything interesting: if $G=K_n$, where $f_G(t)=(1+t)^n$ and 
$e_G(t) = \log(f_G(t) = n \log(1+t)$. On the sub-ring $\mathbb{C}$ of $\mathcal{C}$
generated by complete graphs, we would have then to have $e_G(t) = r(G) \log(1+t)$ with the understanding
that $r(G)$ is the complex number associated to the graph $G$ generated by $1=K_1$. 
We have then $f_G(t)=(1+t)^{r(G)}$ and $\chi(G)=1-f_G(-1) = 1$.
So, on the ``scalar" field $\mathbb{Q}$ or even its completion $\mathbb{R}$ 
generated by complete graphs, the Euler characteristic is always $1$. This is not surprising 
as $\chi(G)=1$ for any complete graph. 

\paragraph{}
While this fuels hope, things fail in general. We can not define Euler 
characteristic consistently in $\mathcal{Q}$. The reason is that
there is no compatibility with the multiplication $*$ in general. The following example 
illustrates this: 

\paragraph{}
Given $G=S^0$, the zero dimensional sphere with two vertices and no edges, we have 
$f_{S^0}(t)=(1+2t)$ and $\chi(S^0)=1-(-1) = 2$, we have 
$\log(f_{S^n})=n \log(1+2t)$ and $\chi(S^n) = 1-e^{n i \pi} = 1+(-1)^n$. So far so good. 
As for multiplication, we have $S^0*S^0 = P^4$ which is the edge-less graph with 4 vertices
 which has Euler characteristic $\chi(G)=4$. We get $e_G(t)=\log(f_{G}(t)) = (\log(1+2t))^2$ and 
$\chi(G) = 1-e^{e_G(-1)} = 1-e^{-\pi^2}$. This is not equal to $4$, 
which is the Euler characteristic of the product $C_4 = S_0 * S_0$. 

\section{Remarks}

\paragraph{}
We originally started to look into this topic by probing more arithmetic-dimension compatibility. 
The clique number $c(G)$ is associated to the maximal dimension as it is the maximal dimension plus 1. While the 
{\bf augmented inductive dimension} ${\rm dim}^+=1+{\rm dim}$ with inductive dimension ${\rm dim}$ 
works with addition by a recent result of Betre and Salinger, the corresponding norm is not compatible 
with graphs. The exact relation of inductive dimension with arithmetic is not yet much explored. It is not excluded that there are 
other constructs, in which algebra, dimension and Euler characteristic interplay in some interesting way. 

\paragraph{}
While the dual ring to $(\mathcal{Z},+,\star)$ is algebraically isomorphic to $\mathcal{Z}$,
the clique number $c$ is not a translational invariant quantity even on that dual ring
because $c(A \cup B) = {\rm max} (c(A),c(B))$. But in principle, one can also use the dual
{\bf ``pebble picture"}, where the point graphs (graphs without edges) play the role of the
classical natural numbers. The fact that the join operation preserves spheres
and the class Dehn-Sommerville graphs tilts towards using the ring structure in which
the join is the addition. The compatibility of Euler characteristic and cohomology makes
the dual story attractive, where the addition is the disjoint union and the multiplication 
is the strong multiplication.

\paragraph{}
There is no way to make $\mathcal{R}$ into a {\bf division algebra}, as the structure of
commutative associative division algebras is known to be either $\mathbb{R}$ or $\mathbb{C}$ by Mazur's theorem.
If every non-zero element was invertible, then the algebra would have to be $\mathbb{R}$ or $\mathbb{C}$
as these are only two real Banach division algebras.
As mentioned in \cite{ArithmeticGraphs} one could look at the smallest field containing $\mathcal{Z}$.
We prefer to keep it a ring as we want to complete it and because in the completion we can not have a field: 
no infinite dimensional Banach algebra can be a field; it would have to be an associative division algebra
over the reals and so $\mathbb{R}$ or $\mathbb{C}$. 
(If it was a field, then $|b|=|a*(b/a)| \leq |a| |b/a| \leq |a| |b| |1/a|$ would imply $|a| |1/a|=1$ and so $|1/a|=1/|a|$ leading
to an equality in all cases implying $|a*c| = |a| |c|$ for all $a,c$ which is the division algebra property.)

\paragraph{}
Can the Banach algebra $\mathcal{R}$ be made into a Hilbert space? 
One can try to define an {\bf inner product} as follows: given two graphs $G,H \in \mathcal{N}$, 
define $\langle G,H \rangle = c(G) c(H)$. 
The additive compatibility of $c$ with the join addition $+$ and the multiplictive compatibility of $c$ with 
scalar multiplication $G \to \lambda * G$, (where $\lambda$ is an element in the completion of the ring generated
by complete graphs which serves as the field over which the vector field is defined), 
we have $\langle \lambda G,H \rangle = \lambda \langle G,H \rangle$ and 
$\langle G+H,K \rangle = \langle G,K \rangle + \langle H,K \rangle$ (and similarly for the right hand side) 
establishing the bilinear property for the inner product. 

\paragraph{}
Dwelling more on the question of a Hilbert space, the ``norm" defined by this inner product
can be zero also on objects which are not zero:
Define $||G|| = \sqrt{ c(G) c(G) } = |c(G)|$. There are many non-zero graphs for which $c(G)=0$, as defined.
For example, for $c(C_4 - C_5)=c(C_4)-c(C_5)=2-2=0$ with cyclic graphs $C_k$, we have $|C_4-C_5| = c(C_4) + c(C_5)=4$
as there are no common additive ``factors" but $||C_4-C_5|| = |c(C_4)-c(C_5)| = |2-2|=0$.
In other words, the bilinear operation $\langle G,H \rangle$ only defines a semi-norm. This is not the end of the
rope as one can then take the {\bf Kolmogorov quotient} and look at equivalence classes identifying
objects $A,B$ with $||A-B||=0$ and get so a Hilbert space. We do not look at this yet here but it might be something
to consider when trying to run a dynamical systems like some quantum mechanics on the space.

\paragraph{}
We can try to do more calculus in the larger playground $\mathcal{Q}$ extending $\mathbb{Q}$.
Here is a simple example, where we look at a dynamical systems like the Fibonacci sequence on graphs.
The graph $G_n$ is recursively is defined as $G_{n+1} = G_n + G_{n-1}$. 
If we start with $G_0=G_1 = S^0$, the $0$-sphere, then $G_n$ are all spheres. The quotients
$G_n/G_{n-1}$ converge to a ``golden mean" in the Banach algebra. 

\begin{figure}[!htpb]
\scalebox{1.1  }{\includegraphics{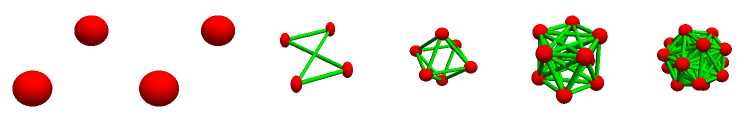}}
\scalebox{1.1  }{\includegraphics{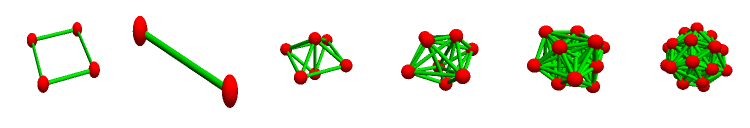}}
\label{fibonacci}
\caption{
Two Fibonacci sequences. Pairs of consecutive elements are added to get the next. 
The first example converges to the graph analogue of the golden mean $1/(1+\sqrt{5})/2$. 
The second is an other object as the circular graph is not generated by complete graphs. 
}
\end{figure}

\paragraph{}
Graphs of the form $G+i H$, where $G,H$ are in $\mathcal{Z}$ build the analog of {\bf Gaussian integers} but again
we don't have a unique multiplicative prime factorization, unless we restrict to connected graphs. 
In general, for graph theoretical analogues of number theoretical questions, the answers could be
easier. Like for Goldbach type questions. 
If $G$ is a complete graph, then $G=A*B$ or $G=A+B$ is only possible if both $A,B$ are complete graphs.
One can now ask what the possible set of graphs $G+H$ is, if $G,H$ are multiplicative primes. While we know
everything about additive primes (the graph complement is connected) we know very little about multiplicative 
primes. One can have the impression that there are many and that most graphs are multiplicative primes. But this
is unexplored. 

\paragraph{}
In the book ``numbers and games" \cite{numbersgames},
numbers and ordinals are constructed in a {\bf Dedekind cut} type manner, leading to 
{\bf surreal numbers} and ``games", a larger class of numbers. 
The story of ``surreal numbers" has been popularized in the short novel \cite{KnuthSurreal}. 
There is no relation between ``graphs" and ``games" but there are formal similarities:
to compare, let us in both case write $G=(G^L|G^R)$ for a ``number". In the game
case, this is a recursive set-up, where $G^L,G^R$ are sets of numbers. 
In the graph case $G=(G^L|G^R)$ has on the left the set of vertices of a graph and 
on the right a set of edges of the graph. 

\paragraph{}
The definition of ``games" starts by defining
$0=(\emptyset|\emptyset)$ is similar as to define the zero graph as the zero $0=(\emptyset|\emptyset)$. 
Now, $1 = (\{0\}|\emptyset)$ is the surreal number $1$. In the graph case,
this looks the same, if we think about $0$ as a ``vertex". But then things start to differ.
In the surreal case $2=(\{1\}|\emptyset)$, while in the graph case 
$2=\{ \{1,2\} |, \{ (1,2) \}$ is the complete graph with two elements. As in integer arithmetic,
the number $-2$ has then to be postulated as something which satisfies $-2+2=0$. 

\paragraph{}
The addition on games $G=(G^L|G^R)$ is defined inductively as $G+H=(G^L+H,G+H^L|G^R+H,G+H^R)$.
In the case of graphs $G=(V,E)=(G^L,G^R)$, the addition is defined as 
$G+H=(G^L \cup H^L | G^R \cup G^R \cup E(G^L,H^L) \}$,
where $E$ is an operation giving from two vertex sets an edge set. In both cases, one has a 
monoid now. There is a difference already however. The class of finite simple graphs does not have
to be constructed recursively. It is already given: take a finite set $V$ and a finite set of unordered
pairs in $V \times V$ and get a graph. 

\paragraph{}
The definition of {\bf negative games} $-G=(-G^R|-G^L)$ is intuitive, if one thinks about the left and right
parts as sets of numbers. In the case of ``graphs", the ``negative graphs" are defined more abstractly
using the Grothendieck completion of the monoid operation $+$. For standard arithmetic in $\mathbb{Q}$,
it took a while to get used to this: the father of algebra, the Greek mathematician Diophantus considered 
negative number an ``absurdity". They only started to appear in the 7'th century, in the work of Brahmagupta, 
but they were considered with suspicion until much later, especially in European mathematics. 
Fibonacci in the beginning of the 13'th century allowed negative numbers as debits or losses. 
This is how the topic is motivated still today when taught to school children.
Still it was not uncommon until the 18th century, that negative numbers were considered nonsensical or 
meaningless. A similar hurdle had to be taken when introducing the complex numbers. The hurdle to get to 
real numbers was taken earlier by the Greeks due to the natural appearance of irrational numbers in geometry
as diagonals in a square or circumference of a circle. 

\paragraph{}
A precise set-theoretical construction of negative numbers is based on the process of {\bf group completion}.
The general insight came surprisingly late and was only formulated in an ultimate way with Grothendieck's concept:
given a commutative monoid $M$, one looks at pairs $(x,y) \in M \times M$ and thinks about it 
as $x-y$. Now define the addition $(x,y) + (a,b) = (x+a,y+b)$, which only uses the monoid addition on 
both ``positive" and ``negative" side. Looking at all pairs double-counts numbers because say
$(5,3)=(6,4)$ as $5-3=6-4$. Mathematically one deals with this by imposing an equivalence relation
$(x,y) \sim (a,b)$ if $x+b = y+a$, which again only uses the monoid structure. Exactly the same 
is done when defining ``negative graphs". It is convenient to write a general graph in the additive
group $\mathcal{N}$ of graphs as $G-H$, where $G,H$ are both graphs. But there are pairs of graphs
$X-Y$ and $A-B$ which are considered equivalent. This happens if $X+B=A+Y$. 

\paragraph{}
In the multiplicative monoid,  we first have to extend the multiplication to 
the additive group by defining $(A-B)*(C-D) = (A*C+B*D) - (A*D+B*D)$ then form the ``rational numbers"
by again using pairs of numbers. A graph in $\mathcal{Q}$ now can be written as $(A-B)/(C-D)$, where
$A,B,C,D$ are finite simple graphs for which $c(C-D)=c(C) - c(D)$ is non-zero. 
The addition and multiplication are done as in school arithmetic. 
In the completion, one has to worry about the case $(1/G) (1/H)$, where $G*H=0$ with non-zero $G$ and $H$.
This is not possible if $G,H$ are in $\mathcal{Q}$ as there, we have assumed for an element $G=P/Q$ that $c(Q) \neq 0$. 
The clique number $c$ can not be extended to all fractions $G/H$ as $c(H)$ can be zero on $\mathcal{Z}$. 
One can still define a norm $|G|=|A-B| = c(A)+c(B)$ if $A-B$ is reduced and have that $|G|>0$ for $G \neq 0$. 
This norm can then be extended to {\bf all non-zero rationals} by $|A/B| = |A|/|B|$. 
The product property $c(A*B) = c(A)*c(B)$ for $A$ generated by complete graphs is pivotal for establishing that the
scalar multiplication works as the ``real numbers" $\mathbb{R}$ sit inside the Banach algebra 
$\mathcal{R}$ as defined. 

\paragraph{}
Implementing an algebra structure on a geometry or on spaces of geometries is an old theme in mathematics.
If we think about ``pebbles" as graphs without edges, then adding pebbles is an operation on graphs
already. But this is different from the graph addition as the addition $1+1$ is not the disjoint graph of 
two pebbles, but a complete graph with two vertices. The arithmetic of complete graphs is via graph 
complement equivalent to the pebble addition of our 10 thousand year old ancestors. 

\paragraph{}
The task to build arithmetic for Euclidean space $\mathbb{R}^n$ led to the concept of 
{\bf division algebras} like $\mathbb{C}$ or $\mathbb{H}$ 
if multiplicative compatibility is required. By theorems of Hurwitz, Frobenius and Mazur, there are very few.
Only two commutative $\mathbb{R}$ and $\mathbb{C}$ and one non-commutative $\mathbb{H}$
if associativity is required. This is related to the fact that only in dimensions $n=1,n=2$ and $n=4$
the unit spheres are Lie groups. The case $n=2$ is the only commutative and
$n=4$, the only non-commutative one. It is no surprise that nature likes these groups when
implementing fundamental forces. When taking the spaces $\mathbb{R}^n$ themselves as elements,
then one has the Euclidean product as well as the tensor product.

\paragraph{}
Ring structures on geometries have been implemented elsewhere in geometry. One can take tensor products
vector bundles on manifolds. Having a division algebra on finite dimensional vector space is related to the
existence whether spheres are Lie groups. In general, the concept of Lie group is a marriage between
algebra and geometry. In algebraic geometry, the concept of {\bf divisors} imposes a group
structure on geometry if one sees the divisor as a geometric object on which
quantities like Euler characteristic are evaluated $\chi(G)+{\rm deg}(D)$. 
In analysis, {\bf differential complexes} are used like de de Rham complex. 
The index plays there the role of the Euler characteristic.
In any case, concepts like divisor algebras, Lie theory, vector bundles, differential complexes divisors 
have all shown that doing algebra on a geometry can be useful. It is even in the name of algebraic geometry
or algebraic topology. Much of commutative algebra is ring theory. 

\paragraph{}
One could also look at different base fields. The fact that there are no division algebras relies
on the chosen base field. One can look at other fields by identifications of complete
graphs. For example, if the complete graph with $3$ elements is identified with the zero graph, 
then one obtains a graph arithmetic over the field $\mathbb{Z}_3$. We have not looked whether one could
get like that to division algebras.

\vfill

\pagebreak

\section{Code}

\paragraph{}
Here is some Mathematica code (check the ArXiv to copy paste): it is not the most efficient way to 
implement things but the slow procedures indicate what is happening better.
It should be readable pseudo code, but it runs and allows to experiment with the basic arithmetic. 

\begin{tiny}
\lstset{language=Mathematica} \lstset{frameround=fttt}
\begin{lstlisting}[frame=single]
NormalizeGraph[s_]:=Module[{e,v,e2,vl,nn,r,V1=VertexList[s],e1=EdgeRules[s],ss},
 e2={};Do[If[Not[e1[[k,1]]==e1[[k,2]]],e2=Append[e2,e1[[k]]]],{k,Length[e1]}]; e1=e2;
 vl=VertexList[s]; nn=Length[vl];r=Table[vl[[k]]->k,{k,nn}];e=e1 /. r; v=V1 /. r;
 UndirectedGraph[Graph[v,e]]];
CliqueNumber[s_]:=Length[First[FindClique[s]]];

ZykovJoin[ss1_,ss2_] := Module[{s1,s2,v1, v2, n1, n2, e1, e2,v,e},
   s1 = NormalizeGraph[ss1]; s2 = NormalizeGraph[ss2];
   v1 = VertexList[s1];      n1 = Length[v1]; v2 = VertexList[s2]+n1;n2=Length[v2];
   e1 = EdgeList[s1];        e2 = EdgeList[s2];
   e1 = If[Length[e1]==0,{},Table[e1[[k,1]] -> e1[[k, 2]], {k, Length[e1]}]];
   e2 = If[Length[e2]==0,{},Table[e2[[k,1]]+ n1->e2[[k,2]]+n1,{k,Length[e2]}]];
   v = Union[v1, v2]; e = Flatten[ Union[{e1, e2, Flatten[Table[
        v1[[k]] -> v2[[l]], {k, Length[v1]}, {l, Length[v2]}]]}]];
   NormalizeGraph[UndirectedGraph[Graph[v, e]]]];
ZykovJoin[ss1_,ss2_,ss3_]:=ZykovJoin[ss1,ZykovJoin[ss2,ss3]];

ZykovProduct[s1_,s2_]:=Module[{v1,v2,e1,e2,v,e={}},
  v1 = VertexList[s1];        v2 = VertexList[s2];
  e1 = Union[EdgeList[s1]];   e2 = Union[EdgeList[s2]];
  e1 = Table[Sort[{e1[[k,1]], e1[[k,2]]}], {k, Length[e1]}];
  e2 = Table[Sort[{e2[[k,1]], e2[[k,2]]}], {k, Length[e2]}];
  v = Partition[Flatten[Table[{v1[[k]],v2[[l]]},{k,Length[v1]},{l,Length[v2]}]],2];
  Do[Do[Do[e = Append[e, {e1[[k, 1]], v2[[m]]} -> {e1[[k, 2]], v2[[l]]}],
      {k,Length[e1]}], {m,Length[v2]}], {l,Length[v2]}];
  Do[Do[Do[e = Append[e, {v1[[m]], e2[[k,1]]} -> {v1[[l]], e2[[k,2]]}],
      {k,Length[e2]}], {m,Length[v1]}], {l,Length[v1]}];
  NormalizeGraph[UndirectedGraph[Graph[v,e]]] ];
ZykovProduct[s1_,s2_,s3_]:=ZykovProduct[s1,ZykovProduct[s2,s3]]; GC=GraphComplement;
StrongProduct[s1_,s2_]:=NormalizeGraph[GC[ZykovProduct[GC[s1],GC[s2]]]];
AdditivePrimeQ[s_]:=ConnectedGraphQ[GraphComplement[s]]; 

ErdoesRenyi[M_,p_]:=Module[{q,e,a},V=Range[M];e=EdgeRules[CompleteGraph[M]]; q={};
  Do[If[Random[]<p,q=Append[q,e[[j]]]],{j,Length[e]}];UndirectedGraph[Graph[V,q]]];

s1=ErdoesRenyi[8,0.4]; s2=ErdoesRenyi[9,0.5]; s3=ErdoesRenyi[10,0.4];
s23=ZykovJoin[s2,s3]; p13=ZykovProduct[s1,s3]; p12=ZykovProduct[s1,s2];
A=ZykovProduct[s1,s23];    B=ZykovJoin[p12,p13];

Print["Distributivity: ",IsomorphicGraphQ[A,B]];
Print["Multiplicativity: ",CliqueNumber[p12]==CliqueNumber[s1]*CliqueNumber[s2]];
Print["Additivity: ",CliqueNumber[s23]==CliqueNumber[s2]+CliqueNumber[s3]];
Print["Additive primes: ",Table[AdditivePrimeQ[CycleGraph[k]], {k, 3, 7}]];
\end{lstlisting}
\end{tiny}

\begin{figure}[!htpb]
\scalebox{0.516}{\includegraphics{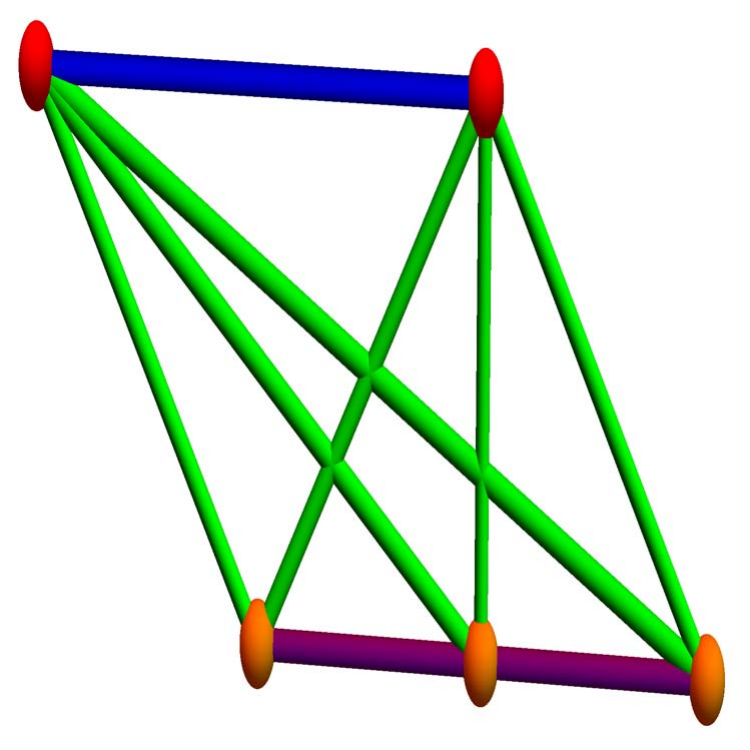}}
\scalebox{0.516}{\includegraphics{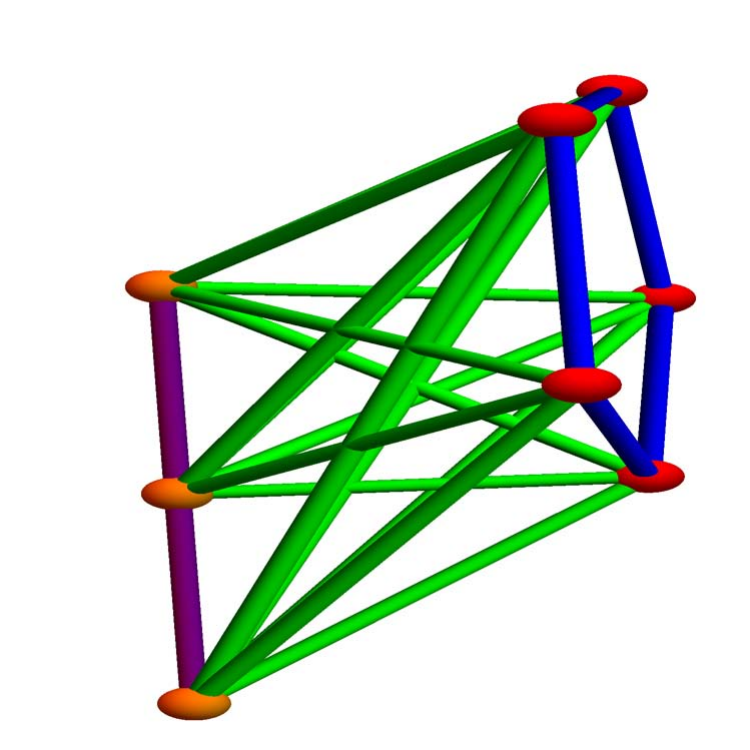}}
\scalebox{0.516}{\includegraphics{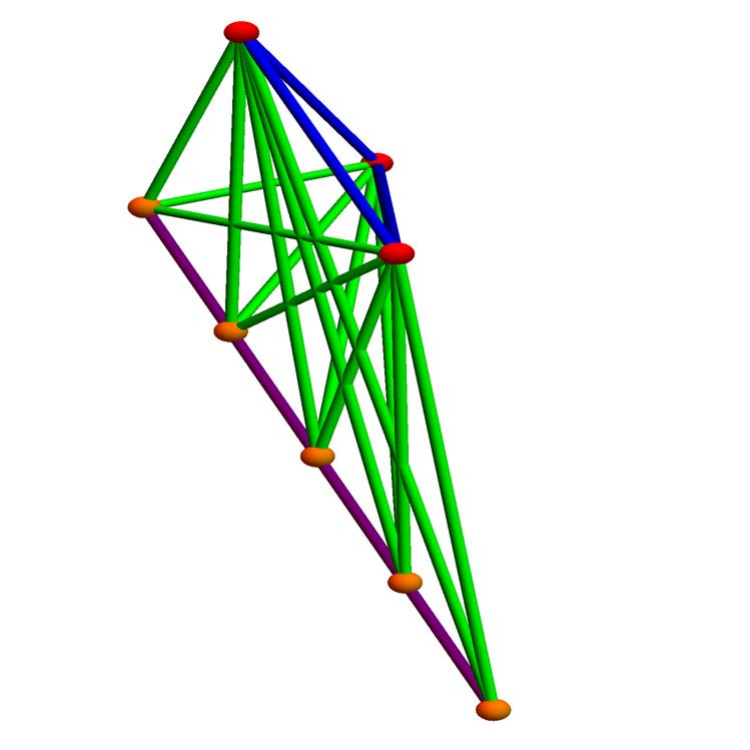}}
\scalebox{0.516}{\includegraphics{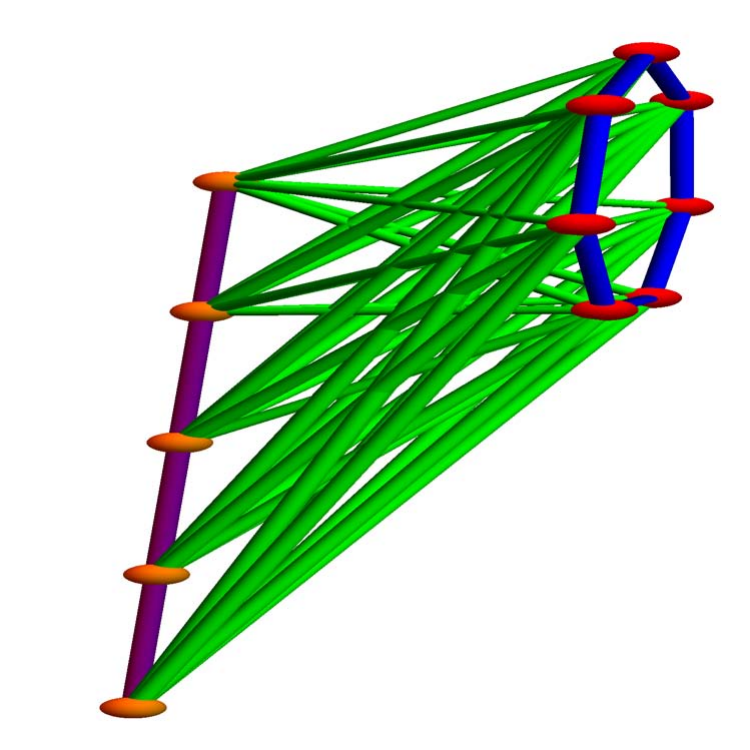}}
\scalebox{0.516}{\includegraphics{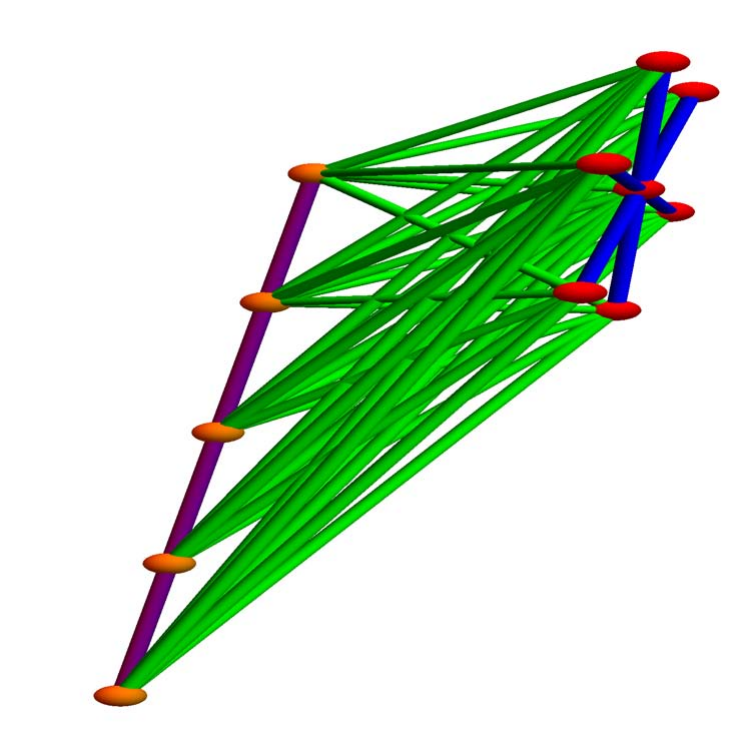}}
\scalebox{0.516}{\includegraphics{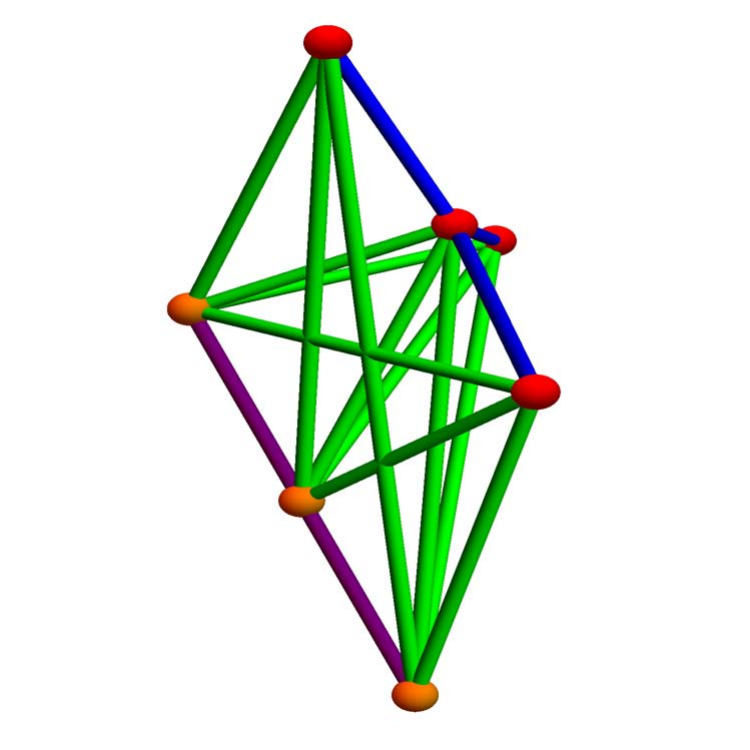}}
\label{addition}
\caption{
The join addition $G+H$ is obtained by taking the disjoint union of the graphs 
followed by joining all mixed vertex pairs. 
}
\end{figure}

\begin{figure}[!htpb]
\scalebox{0.516}{\includegraphics{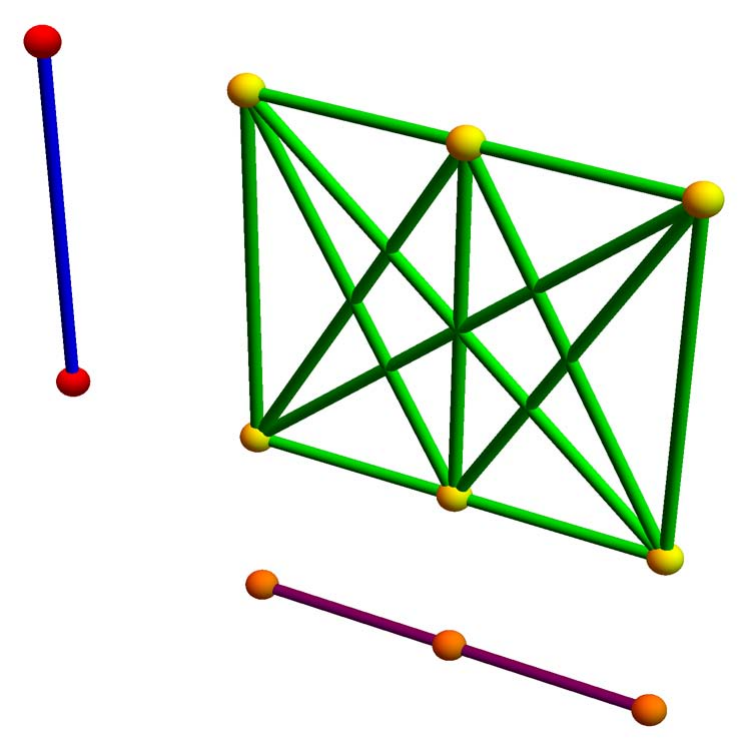}}
\scalebox{0.516}{\includegraphics{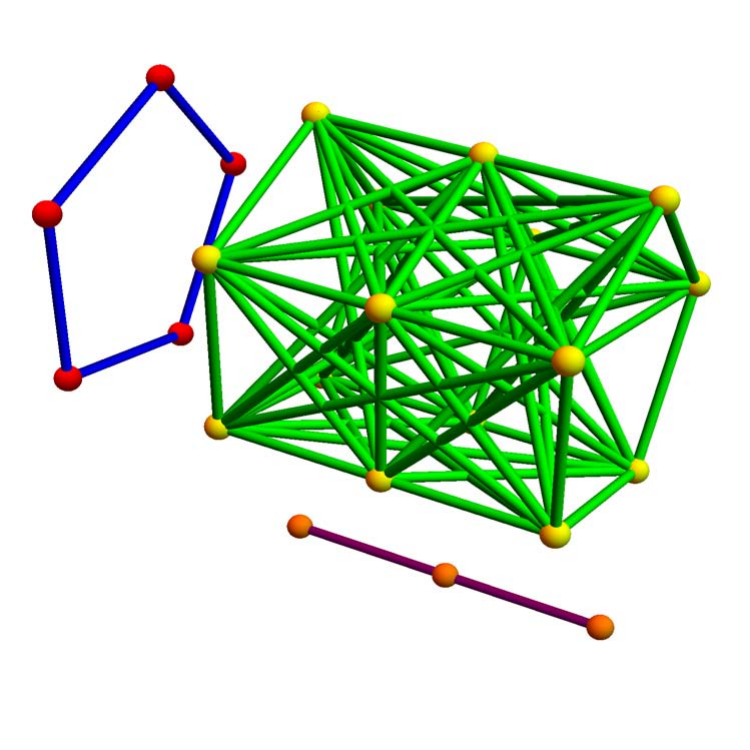}}
\scalebox{0.516}{\includegraphics{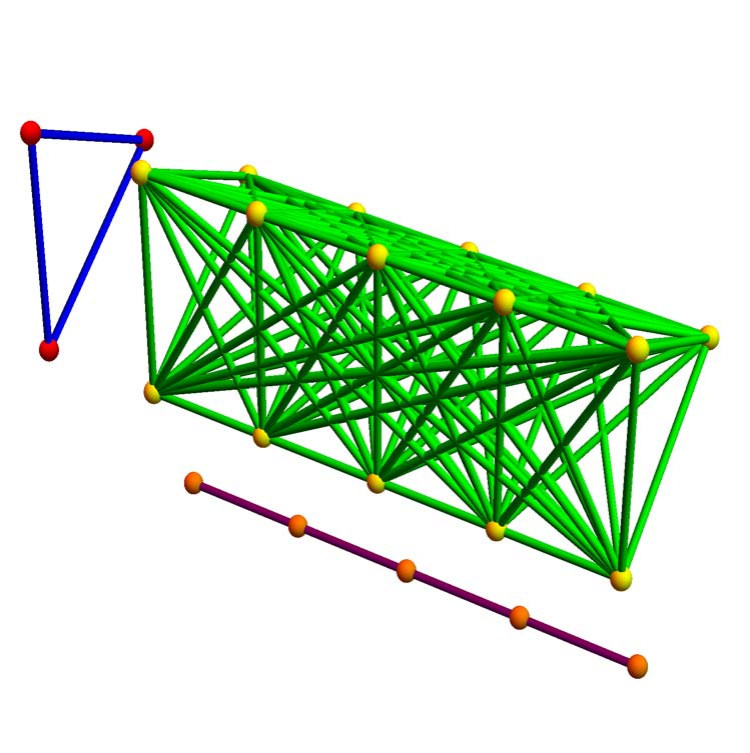}}
\scalebox{0.516}{\includegraphics{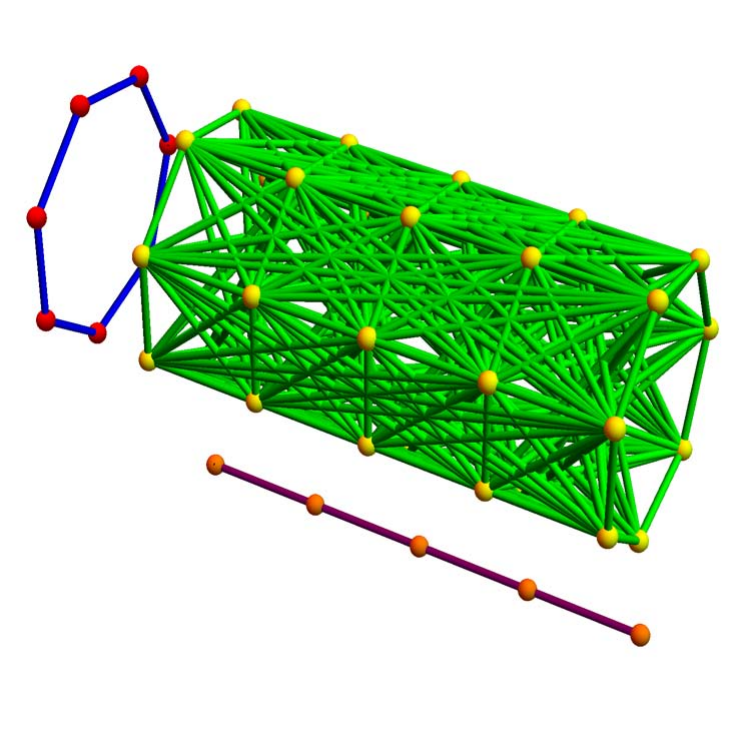}}
\scalebox{0.516}{\includegraphics{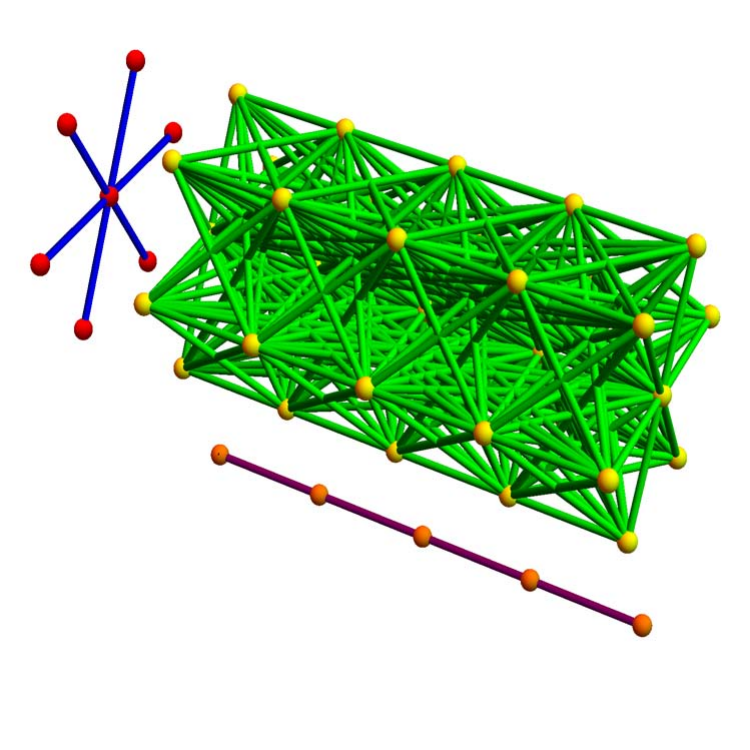}}
\scalebox{0.516}{\includegraphics{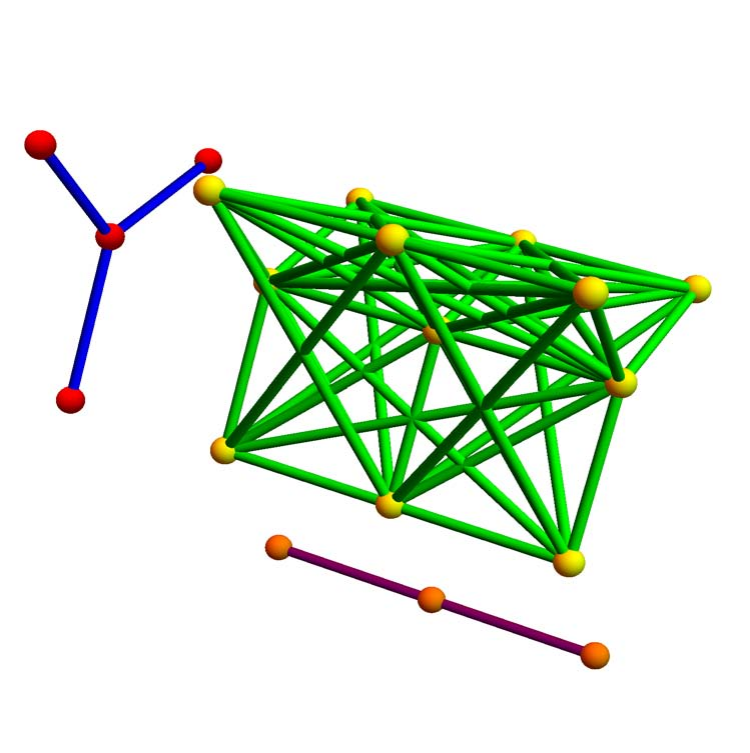}}
\label{multiplication}
\caption{
The multiplication $G *H$ has the product set $V \times W$ as vertex sets $V$ and $W$
of $G$ and $H$ and connects two vertices with an edge, 
if the projection to one of the two graphs is an edge there. 
}
\end{figure}

\begin{figure}[!htpb]
\scalebox{0.2}{\includegraphics{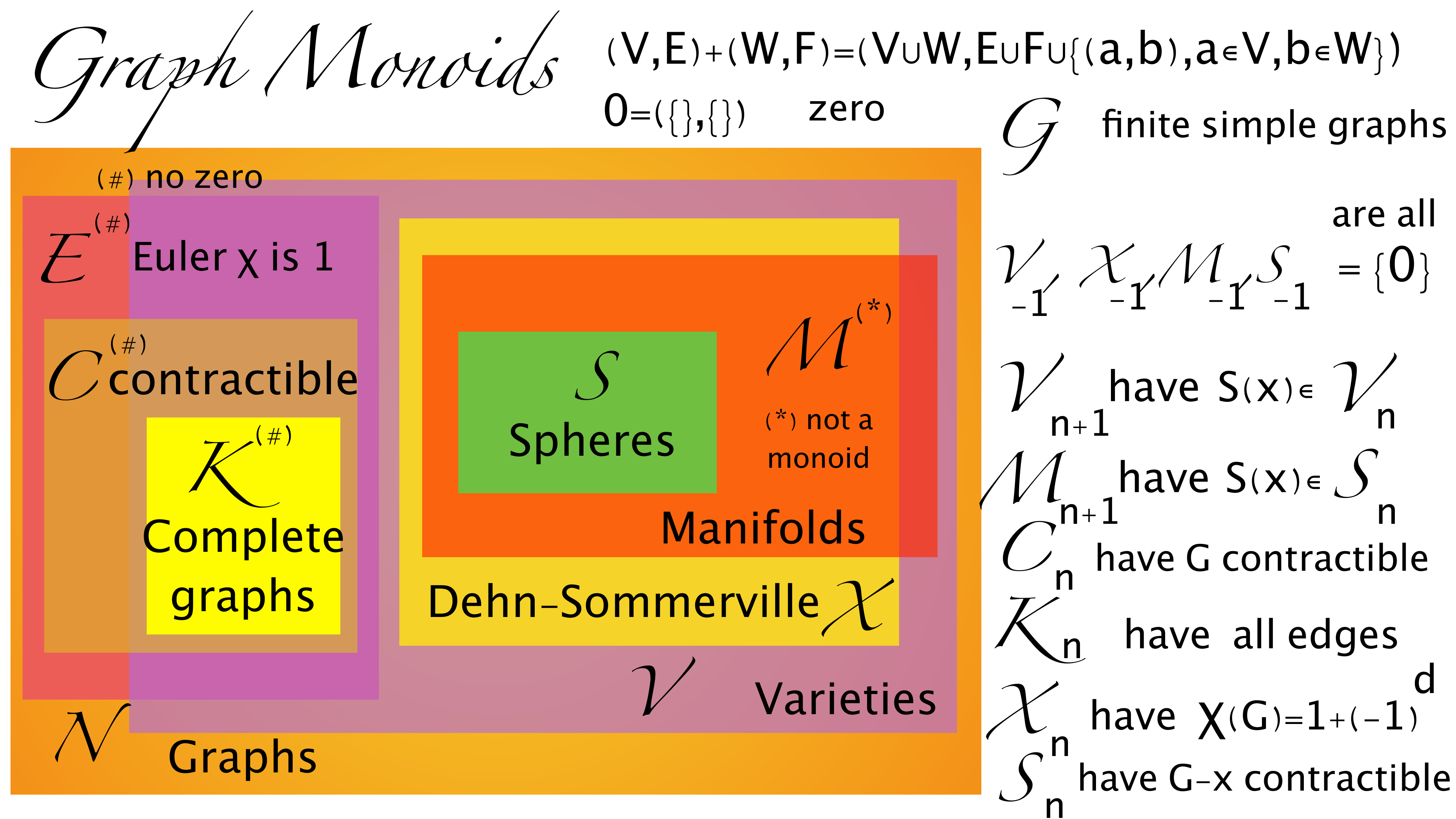}}
\label{monoids}
\caption{
Some Monoids of graphs. The spheres, the Dehn-Sommerville graphs, 
the varieties and all the finite simple graphs are all monoids. 
We call $\mathcal{N}$ the set of all finite simple graphs and extend it to 
$\mathcal{Z}, \mathcal{Q}, \mathcal{R}$ which is then a Banach algebra. 
}
\end{figure}

\bibliographystyle{plain}

\end{document}